\documentclass[12pt,leqno]{article}
\usepackage{amssymb}
\usepackage[mathscr]{eucal}
\usepackage{amsmath,amssymb,latexsym,theorem,bbm}
\usepackage{color}
\usepackage{appendix}

\newcommand{\DS}{\displaystyle}

\newcommand{\SC}{\scriptstyle}

\newcommand{\CC}{\mathbb{C}}

\newcommand{\RR}{\mathbb{R}}

\newcommand{\ta}{\widetilde{a}}

\newcommand{\bh}{{\boldsymbol{h}}}
\newcommand{\bI}{{\boldsymbol{I}}}
\newcommand{\bJ}{{\boldsymbol{J}}}

\newcommand{\br}{{\boldsymbol{r}}}

\newcommand{\bDelta}{{\boldsymbol{\Delta}}}

\newcommand{\btheta}{{\boldsymbol{\theta}}}

\newcommand{\bzero}{{\boldsymbol{0}}}

\newcommand{\cA}{{\mathcal A}}
\newcommand{\cB}{{\mathcal B}}

\newcommand{\cD}{{\mathcal D}}
\newcommand{\cE}{{\mathcal E}}

\newcommand{\cN}{{\mathcal N}}

\newcommand{\cY}{{\mathcal Y}}
\newcommand{\cZ}{{\mathcal Z}}
\newcommand{\cW}{{\mathcal W}}

\newcommand{\dd}{\mathrm{d}}
\newcommand{\ee}{\mathrm{e}}

\newcommand{\EE}{\operatorname{\mathbb{E}}}
\newcommand{\PP}{{\operatorname{\mathbb{P}}}}

\newcommand{\oo}{\operatorname{o}}
\newcommand{\OO}{\operatorname{O}}

\renewcommand{\Re}{\operatorname{Re}}
\renewcommand{\Im}{\operatorname{Im}}

\newcommand{\ha}{\widehat{a}}

\renewcommand{\leq}{\leqslant}
\renewcommand{\geq}{\geqslant}

\newcommand{\stoch}{\stackrel{\PP}{\longrightarrow}}
\newcommand{\distr}{\stackrel{\cD}{\longrightarrow}}

\newcommand{\distre}{\stackrel{\cD}{=}}

\newcommand{\as}{\stackrel{{\mathrm{a.s.}}}{\longrightarrow}}

\newcommand{\proofend}{\hfill\mbox{$\Box$}}

\numberwithin{equation}{section}

\theoremstyle{change} \theorembodyfont{\em}
\newtheorem{Lem}{Lemma.}[section]

\newtheorem{Pro}[Lem]{Proposition.}
\newtheorem{Cor}[Lem]{Corollary.}
\newtheorem{Def}[Lem]{Definition.}

\theorembodyfont{\rm}
\newtheorem{Rem}[Lem]{Remark.}

\begin{document}

\begin{center}
 {\bfseries\Large Asymptotic inference for a stochastic differential} \\[2mm]
 {\bfseries\Large equation with uniformly distributed time delay} \\[5mm]
 {\sc\large J\'anos Marcell $\text{Benke}^*$ \ and \ Gyula Pap}
\end{center}

\vskip0.2cm

\noindent
 Bolyai Institute, University of Szeged,
 Aradi v\'ertan\'uk tere 1, H--6720 Szeged, Hungary.

\noindent e--mails: jbenke@math.u-szeged.hu (J. M. Benke),
                    papgy@math.u-szeged.hu (G. Pap).

\noindent * Corresponding author.

\renewcommand{\thefootnote}{}

\vspace*{-7mm}


\begin{abstract}
For the affine stochastic delay differential equation
 \[
   \dd X(t) = a \int_{-1}^0 X(t + u) \, \dd u \, \dd t + \dd W(t) ,
   \qquad t \geq 0 ,
 \]
 the local asymptotic properties of the likelihood function are studied.
Local asymptotic normality is proved in case of
 \ $a \in \bigl(-\frac{\pi^2}{2}, 0\bigr)$, \ local asymptotic mixed normality
 is shown if \ $a \in \bigl(0, \infty\bigr)$, \ periodic local asymptotic
 mixed normality is valid if \ $a \in \bigl(-\infty, -\frac{\pi^2}{2}\bigr)$,
 \ and only local asymptotic quadraticity holds at the points \ $-\frac{\pi^2}{2}$
 \ and \ $0$. 
\ Applications to the asymptotic behaviour of the maximum likelihood estimator
 \ $\ha_T$ of \ $a$ \ based on \ $(X(t))_{t\in[0,T]}$ \ are given \ as \ $T \to \infty$.
\end{abstract}

\section{Introduction}

Assume \ $(W(t))_{t\in\RR_+}$ \ is a standard Wiener process, \ $a \in \RR$,
 \ and \ $(X^{(a)}(t))_{t\in\RR_+}$ \ is a solution of the affine stochastic
 delay differential equation (SDDE)
 \begin{equation}\label{SDDE}
  \begin{cases}
   \dd X(t) = a \int_{-1}^0 X(t + u) \, \dd u \, \dd t + \dd W(t) ,
    & t \in \RR_+ , \\
   X(t) = X_0(t) , & t \in [-1, 0] , 
  \end{cases}
 \end{equation}
 where \ $(X_0(t))_{t\in[-1,0]}$ \ is a continuous stochastic process
 independent of \ $(W(t))_{t\in\RR_+}$.
\ The SDDE \eqref{SDDE} can also be written in the integral form
 \begin{equation}\label{iSDDE}
  \begin{cases}
   X(t) = X_0(0) + a \int_0^t \int_{-1}^0 X(s + u) \, \dd u \, \dd s + W(t) ,
    & t \in \RR_+ , \\
   X(t) = X_0(t) , & t \in [-1, 0] . 
  \end{cases}
 \end{equation}
Equation \eqref{SDDE} is a special case of the affine stochastic delay differential
 equation
 \begin{equation}\label{ASDDE}
  \begin{cases}
   \dd X(t) = \int_{-r}^0 X(t + u) \, m_\theta(\dd u) \, \dd t + \dd W(t) ,
    & t \in \RR_+ , \\
   X(t) = X_0(t) , & t \in [-r, 0] , 
  \end{cases}
 \end{equation}
 where \ $r > 0$, \ and \ for each \ $\theta \in \Theta$, \ $m_\theta$, \ is a
 finite signed measure on \ $[-r, 0]$ \ see Gushchin and K\"uchler \cite{GusKuc2003}.
In that paper local asymptotic normality has been proved for stationary solutions.
In Gushchin and K\"uchler \cite{GusKuc1999}, the special case of \eqref{ASDDE} has
 been studied with \ $r = 1$, \ $\Theta = \RR^2$, \ and
 \ $m_\theta = a \delta_0 + b \delta_{-1}$ \ for \ $\theta = (a, b)$, \ where
 \ $\delta_x$ \ denotes the Dirac measure concentrated at \ $x \in \RR$, \ and they
 described the local properties of the likelihood function for the whole
 parameter space \ $\RR^2$.
 
The solution \ $(X^{(a)}(t))_{t\in\RR_+}$ \ of \eqref{SDDE} exists, is
 pathwise uniquely determined and can be represented as
 \begin{equation}\label{solution}
  X^{(a)}(t)
  = x_{0,a}(t) X_0(0)
    + a \int_{-1}^0 \int_u^0 x_{0,a}(t + u - s) X_0(s) \, \dd s \, \dd u
    + \int_0^t x_{0,a}(t - s) \, \dd W(s) ,
 \end{equation}
 for \ $t \in \RR_+$, \ where \ $(x_{0,a}(t))_{t\in[-1,\infty)}$ \ denotes
 the so-called fundamental solution of the deterministic homogeneous delay
 differential equation
 \begin{equation}\label{fund}
  \begin{cases}
   x(t) = x_0(0) + a \int_0^t \int_{-1}^0 x(s + u) \, \dd u \, \dd s ,
    & t \in \RR_+ , \\
   x(t) = x_0(t) , & t \in [-1, 0] . 
  \end{cases}
 \end{equation}
 with initial function
 \[
   x_0(t) := \begin{cases}
              0 , & t \in [-1, 0) , \\
              1 , & t = 0 .
             \end{cases}
 \]
In the trivial case of \ $a = 0$, \ we have \ $x_{0,0}(t) = 1$,
 \ $t \in \RR_+$, \ and \ $X^{(0)}(t) = X_0(0) + W(t)$, \ $t \in \RR_+$.
\ In case of \ $a \in \RR \setminus \{0\}$, \ the behaviour of
 \ $(x_{0,a}(t))_{t\in[-1,\infty)}$ \ is connected with the so-called
 characteristic function \ $h_a : \CC \to \CC$, \ given by
 \begin{equation}\label{char1}
  h_a(\lambda) := \lambda - a \int_{-1}^0 \ee^{\lambda u} \, \dd u ,
  \qquad \lambda \in \CC ,
 \end{equation}
 and the set \ $\Lambda_a$ \ of the (complex) solutions of the 
 so-called characteristic equation for \eqref{fund},
 \begin{equation}\label{char2}
  \lambda - a \int_{-1}^0 \ee^{\lambda u} \, \dd u = 0 .
 \end{equation}
Note that a complex number \ $\lambda$ \ solves \eqref{char2} if and
 only if \ $(\ee^{\lambda t})_{t\in[-1,\infty)}$ \ solves \eqref{fund}
 with initial function \ $x_0(t) = \ee^{\lambda t}$, \ $t \in [-1, 0]$.
\ Applying usual methods (e.g., argument principle in complex analysis
 and the existence of local inverses of holomorphic functions), one
 can derive the following properties of the set \ $\Lambda_a$, \ see,
 e.g., Rei{\ss} \cite{Rei}.
We have \ $\Lambda_a \ne \emptyset$, \ and \ $\Lambda_a$ \ consists of
 isolated points.
Moreover, \ $\Lambda_a$ \ is countably infinite, and for each
 \ $c \in \RR$, \ the set
 \ $\{\lambda \in \Lambda_a : \Re(\lambda) \geq c\}$ \ is finite.
In particular,
 \[
   v_0(a) := \sup\{\Re(\lambda) : \lambda \in \Lambda_a\} < \infty .
 \]
Put
 \[
   v_1(a)
   := \sup\{\Re(\lambda)
            : \lambda \in \Lambda_a , \, \Re(\lambda) < v_0(a)\} ,
 \]
 where \ $\sup \emptyset := - \infty$.
\ We have the following cases:
 \renewcommand{\labelenumi}{{\rm(\roman{enumi})}}
 \begin{enumerate}
  \item 
   If \ $a \in \bigl(-\frac{\pi^2}{2}, 0\bigr)$ \ then \ $v_0(a) < 0$;
  \item 
   If \ $a = -\frac{\pi^2}{2}$ \ then \ $v_0(a) = 0$ \ and
    \ $v_0(a) \notin \Lambda_a$;
  \item 
   If \ $a \in \bigl(-\infty, -\frac{\pi^2}{2}\bigr)$ \ then
    \ $v_0(a) > 0$ \ and \ $v_0(a) \notin \Lambda_a$;
  \item 
   If \ $a \in (0, \infty)$ \ then \ $v_0(a) > 0$,
    \ $v_0(a) \in \Lambda_a$, \ $m(v_0(a)) = 1$ \ (where \ $m(v_0(a))$
    \ denotes the multiplicity of \ $v_0(a)$), and \ $v_1(a) < 0$.
 \end{enumerate}
For any \ $\gamma > v_0(a)$, \ we have
 \ $x_{0,a}(t) = \OO(\ee^{\gamma t})$, \ $t \in \RR_+$.
\ In particular, \ $(x_{0,a}(t))_{t\in\RR_+}$ \ is square integrable
 if (and only if, see Gushchin and K\"uchler \cite{GusKuc2000})
 \ $v_0(a) < 0$. 
\ The Laplace transform of \ $(x_{0,a}(t))_{t\in\RR_+}$ \ is given by
 \[
   \int_0^\infty \ee^{-\lambda t} x_{0,a}(t) \, \dd t
   = \frac{1}{h_a(\lambda)} , \qquad \lambda \in \CC , \qquad
   \Re(\lambda) > v_0(a) .   
 \]
Based on the inverse Laplace transform and Cauchy's residue theorem,
 the following crucial lemma can be shown (see, e.g., Gushchin and
 K\"uchler \cite[Lemma 1.1]{GusKuc1999}). 
 
\begin{Lem}\label{res}
For each \ $a \in \RR \setminus \{0\}$ \ and each
 \ $c \in (-\infty, v_0(a))$, \ there exists \ $\gamma \in (-\infty, c)$
 \ such that the fundamental solution \ $(x_{0,a}(t))_{t\in[-1,\infty)}$
 \ of \eqref{fund} can be represented in the form
 \[
   x_{0,a}(t)
   = \psi_{0,a}(t) \ee^{v_0(a) t}
     + \sum_{\underset{\SC \Re(\lambda)\in[c,v_0(a))}{\lambda \in \Lambda_a}}
        c_a(\lambda) \ee^{\lambda t}
     + \oo(\ee^{\gamma t}) , \qquad \text{as \ $t \to \infty$,}
 \]
 with some constants \ $c_a(\lambda)$, \ $\lambda \in \Lambda_a$, \ and with
 \[
   \psi_{0,a}(t)
   := \begin{cases}
       {\DS \frac{v_0(a)}{v_0(a)^2+2v_0(a)-a}},
        & \text{if \ $v_0(a) \in \Lambda_a$ \ and \ $m(v_0(a)) = 1$,} \\[3mm]
       A_0(a) \cos(\kappa_0(a) t) + B_0(a) \sin(\kappa_0(a) t)
        & \text{if \ $v_0(a) \notin \Lambda_a$,}
      \end{cases}
 \]
 with \ $\kappa_0(a) := |\Im(\lambda_0(a))|$, \ where
 \ $\lambda_0(a) \in \Lambda_a$ \ is given by
 \ $\Re(\lambda_0(a)) = v_0(a)$, \ and
 \begin{align*}
  A_0(a)
  &:= \frac{2[(v_0(a)^2-\kappa_0(a)^2)(v_0(a)-2)-av_0(a)]}
           {(v_0(a)^2-\kappa_0(a)^2+2v_0(a)-a)^2
            +4\kappa_0(a)^2(v_0(a)+1)^2} , \\
  B_0(a)
  &:= \frac{2(v_0(a)^2+\kappa_0(a)^2+a)\kappa_0(a)}
           {(v_0(a)^2-\kappa_0(a)^2+2v_0(a)-a)^2
            +4\kappa_0(a)^2(v_0(a)+1)^2} .
 \end{align*}
\end{Lem}

\section{Quadratic approximations to likelihood ratios}
\label{mlan}

We recall some definitions and statements concerning quadratic approximations
 to likelihood ratios based on Jeganathan \cite{Jeg}, Le Cam and Yang
 \cite{LeCamYang} and van der Vaart \cite{Vaart}.

Let \ $(\Omega, \cA, \PP)$ \ be a probability space.
Let \ $\Theta \subset \RR^p$ \ be an open set.
For each \ $\btheta \in \Theta$, \ let
 \ $(X^{(\btheta)}(t))_{t\in[-1,\infty)}$ \ be a continuous stochastic
 process on \ $(\Omega, \cA, \PP)$.
\ For each \ $T \in \RR_+$, \ let \ $\PP_{\btheta,T}$ \ be the probability
 measure induced by \ $(X^{(\btheta)}(t))_{t\in[-1,T]}$ \ on the space
 \ $(C([-1, T]), \cB(C([-1, T])))$. 
 
\begin{Def}
The family
 \ $(C([-1, T]), \cB(C([-1, T])), \{\PP_{\btheta,T}
     : \btheta \in \Theta\})_{T\in\RR_{++}}$
 \ of statistical experiments is said to have locally asymptotically quadratic
 (LAQ) likelihood ratios at \ $\btheta \in \Theta$ \ if there exist (scaling)
 matrices \ $\br_{\btheta,T} \in \RR^{p\times p}$, \ $T \in \RR_{++}$,
 \ random vectors \ $\bDelta_{\btheta} : \Omega \to \RR^p$ \ and
 \ $\bDelta_{\btheta,T} : \Omega \to \RR^p$, \ $T \in \RR_{++}$, \ and random
 matrices \ $\bJ_{\btheta} : \Omega \to \RR^{p \times p}$ \ and
 \ $\bJ_{\btheta,T} : \Omega \to \RR^{p \times p}$, \ $T \in \RR_{++}$, \ such
 that
 \begin{equation}\label{defLAQ}
  \log \frac{\dd\PP_{\btheta+\br_{\btheta,T}\bh_T,T}}{\dd\PP_{\btheta,T}}
   (X^{(\btheta)}|_{[-1, T]})
  = \bh_T^\top \bDelta_{\btheta,T}
    - \frac{1}{2} \bh_T^\top \bJ_{\btheta,T} \bh_T
    + \oo_{\PP}(1)
  \qquad \text{as \ $T \to \infty$}
 \end{equation}
 whenever \ $\bh_T \in \RR^p$, \ $T \in \RR_{++}$, \ is a bounded family
 satisfying \ $\btheta + \br_{\btheta,T} \bh_T \in \Theta$ \ for all
 \ $T \in \RR_{++}$,
 \begin{equation}\label{LAQD}
  (\bDelta_{\btheta,T}, \bJ_{\btheta,T})
  \distr (\bDelta_{\btheta}, \bJ_{\btheta})
  \qquad \text{as \ $T \to \infty$,}
 \end{equation}
 and we have
 \begin{equation}\label{LAQJ}
  \PP\left(\text{$\bJ_{\btheta}$ is symmetric and strictly positive definite}
     \right)
  = 1   
 \end{equation}
 and
 \begin{equation}\label{LAQDJ}
  \EE\left(\exp\left\{\bh^\top \bDelta_{\btheta} 
                      - \frac{1}{2} \bh^\top \bJ_{\btheta} \bh \right\}\right)
  = 1 , \qquad \bh \in \RR^p .
 \end{equation}
\end{Def}

\begin{Def}
A family
 \ $(C([-1, T]), \cB(C([-1, T])), \{\PP_{\btheta,T}
     : \btheta \in \Theta\})_{T\in\RR_{++}}$
 \ of statistical experiments is said to have locally asymptotically mixed
 normal (LAMN) likelihood ratios at \ $\btheta \in \Theta$ \ if it is LAQ at
 \ $\btheta \in \Theta$, \ and the conditional distribution of
 \ $\bDelta_{\btheta}$ \ given \ $\bJ_{\btheta}$ \ is
 \ $\cN_p(\bzero, \bJ_{\btheta})$, \ or, equivalently, there exist a random
 vector \ $\cZ : \Omega \to \RR^p$ \ and a random matrix
 \ $\eta_\btheta : \Omega \to \RR^{p \times p}$, \ such that they are
 independent, \ $\cZ \distre \cN_p(\bzero, \bI_p)$, \ and
 \ $\bDelta_{\btheta} = \eta_\btheta \cZ$, 
 \ $\bJ_{\btheta} = \eta_\btheta \eta_\btheta^\top$.
\end{Def}

\begin{Def}\label{Def_PLAMN}
The family
 \ $(C([-1, T]), \cB(C([-1, T])), \{\PP_{\btheta,T}
     : \btheta \in \Theta\})_{T\in\RR_{++}}$
 \ of statistical experiments is said to have periodic locally asymptotically
 mixed normal (PLAMN) likelihood ratios at \ $\btheta \in \Theta$ \ if there
 exist \ $D \in \RR_{++}$, \ (scaling) matrices
 \ $\br_{\btheta,T} \in \RR^{p\times p}$, \ $T \in \RR_{++}$, \ random vectors
 \ $\bDelta_{\btheta}(d) : \Omega \to \RR^p$, \ $d \in [0, D)$, \ and
 \ $\bDelta_{\btheta,T} : \Omega \to \RR^p$, \ $T \in \RR_{++}$, \ and
 random matrices \ $\bJ_{\btheta}(d) : \Omega \to \RR^{p \times p}$,
 \ $d \in [0, D)$, \ and \ $\bJ_{\btheta,T} : \Omega \to \RR^{p \times p}$,
 \ $T \in \RR_{++}$, \ such that \eqref{defLAQ} holds whenever
 \ $\bh_T \in \RR^p$, \ $T \in \RR_{++}$, \ is a bounded family satisfying
 \ $\btheta + \br_{\btheta,T} \bh_T \in \Theta$ \ for all \ $T \in \RR_{++}$,
 \begin{equation}\label{LAQDg}
  (\bDelta_{\btheta,kD+d}, \bJ_{\btheta,kD+d})
  \distr (\bDelta_{\btheta}(d), \bJ_{\btheta}(d))
  \qquad \text{as \ $k \to \infty$}
 \end{equation}
 for all \ $d \in [0, D)$, \ we have
 \begin{equation}\label{LAQJg}
  \PP\left(\text{$\bJ_{\btheta}(d)$ is symmetric and strictly positive
                 definite}\right)
  = 1 , \qquad d \in [0, D) ,  
 \end{equation}
 and for each \ $d \in [0, D)$, \ the conditional distribution of
 \ $\bDelta_{\btheta}(d)$ \ given \ $\bJ_{\btheta}(d)$ \ is
 \ $\cN_p(\bzero, \bJ_{\btheta}(d))$, \ or, equivalently, there exist a
 random vector \ $\cZ : \Omega \to \RR^p$ \ and a random matrix
 \ $\eta_\btheta(d) : \Omega \to \RR^{p \times p}$ \ such that they are
 independent, \ $\cZ \distre \cN_p(\bzero, \bI_p)$, \ and
 \ $\bDelta_{\btheta}(d) = \eta_\btheta(d) \cZ$, 
 \ $\bJ_{\btheta}(d) = \eta_\btheta(d) \eta_\btheta^\top(d)$.
\end{Def}

\begin{Rem}
The notion of LAMN is defined in Le Cam and Yang \cite{LeCamYang} and Jeganathan
 \cite{Jeg} so that PLAMN in the sense of Definiton \ref{Def_PLAMN} is LAMN as well.
\end{Rem}

\begin{Def}
A family 
 \ $(C([-1, T]), \cB(C([-1, T])), \{\PP_{\btheta,T}
     : \btheta \in \Theta\})_{T\in\RR_{++}}$
 \ of statistical experiments is said to have locally asymptotically normal
 (LAN) likelihood ratios at \ $\btheta \in \Theta$ \ if it is LAMN at
 \ $\btheta \in \Theta$, \ and \ $\bJ_\btheta$ \ is deterministic.
\end{Def}

\section{Radon--Nikodym derivatives}
\label{section_RN}

From this section, we will consider the SDDE \eqref{SDDE} with fixed
 continuous initial process \ $(X_0(t))_{t\in[-1,0]}$.
\ Further, for all \ $T \in \RR_{++}$, \ let \ $\PP_{a,T}$ \ be the
 probability measure induced by \ $(X^{(a)}(t))_{t\in[-1,T]}$ \ on
 \ $(C([-1, T]), \cB(C([-1, T])))$. 
\ In order to calculate Radon--Nikodym derivatives
 \ $\frac{\dd \PP_{\ta,T}}{\dd \PP_{a,T}}$ \ for certain
 \ $a, \ta \in \RR$, \ we need the following statement, which
 can be derived from formula (7.139) in Section 7.6.4 of Liptser and Shiryaev
 \cite{LipShiI}.

\begin{Lem}\label{RN}
Let \ $a, \ta \in \RR$.
\ Then for all \ $T \in \RR_{++}$, \ the measures \ $\PP_{a,T}$ \ and
 \ $\PP_{\ta,T}$ \ are absolutely continuous with respect to each other,
 and
 \begin{align*}
  &\log \frac{\dd \PP_{\ta,T}}{\dd \PP_{a,T}}(X^{(a)}|_{[-1,T]}) \\
  &= (\ta - a) \int_0^T \int_{-1}^0 X^{(a)}(t+u) \, \dd u \, \dd X^{(a)}(t)
     - \frac{\ta^2 - a^2}{2}
       \int_0^T \left(\int_{-1}^0 X^{(a)}(t+u) \, \dd u\right)^2 \dd t \\
  &= (\ta - a) \int_0^T \int_{-1}^0 X^{(a)}(t+u) \, \dd u \, \dd W(t)
     - \frac{(\ta - a)^2}{2}
       \int_0^T \left(\int_{-1}^0 X^{(a)}(t+u) \, \dd u\right)^2 \dd t .
 \end{align*}
\end{Lem}

In order to investigate convergence of the family
 \begin{equation}\label{cET}
  (\cE_T)_{T\in\RR_{++}}
  := \big( C(\RR_+), \cB(C(\RR_+)) ,
           \{\PP_{a,T} : a \in \RR\} \big)_{T\in\RR_{++}}
 \end{equation}
 of statistical experiments, we derive the following corollary.

\begin{Cor}\label{RN_Cor}
For each \ $a \in \RR$, \ $T \in \RR_{++}$, \ $r_{a,T} \in \RR$ \ and
 \ $h_T \in \RR$, \ we have
 \[
   \log \frac{\dd \PP_{a+r_{a,T}h_T,T}}{\dd \PP_{a,T}}(X^{(a)}|_{[-1,T]})
   = h_T \Delta_{a,T} - \frac{1}{2} h_T^2 J_{a,T} ,
 \]
 with
 \[
   \Delta_{a,T}
   := r_{a,T} \int_0^T \int_{-1}^0 X^{(a)}(t+u) \, \dd u \, \dd W(t) , \qquad
   J_{a,T}
   := r_{a,T}^2
      \int_0^T \left(\int_{-1}^0 X^{(a)}(t+u) \, \dd u\right)^2 \dd t .
 \]
Consequently, the quadratic approximation \eqref{defLAQ} is valid.
\end{Cor}

\section{Local asymptotics of likelihood ratios}
\label{section_LALR}

\begin{Pro}\label{LAN}
If \ $a \in \bigl(-\frac{\pi^2}{2}, 0\bigr)$ \ then the family
 \ $(\cE_T)_{T\in\RR_{++}}$ \ of statistical experiments given in
 \eqref{cET} is LAN at \ $a$ \ with scaling
 \ $r_{a,T} = \frac{1}{\sqrt{T}}$, \ $T \in \RR_{++}$, \ and with
 \[
   J_a
   = \int_0^\infty
      \left(\int_{-1}^0 x_{0,a}(t + u) \, \dd u\right)^2 \dd t . 
 \]
\end{Pro}

\begin{Pro}\label{LAQ1}
The family \ $(\cE_T)_{T\in\RR_{++}}$ \ of statistical experiments given in
 \eqref{cET} is LAQ at \ $0$ \ with scaling \ $r_{0,T} = \frac{1}{T}$,
 \ $T \in \RR_{++}$, \ and with
 \[
   \Delta_0 = \int_0^1 \cW(t) \, \dd \cW(t) , \qquad
   J_0 = \int_0^1 \cW(t)^2 \, \dd t,
 \]
 where \ $(\cW(t))_{t\in[0,1]}$ \ is a standard Wiener process.
\end{Pro}

\begin{Pro}\label{LAQ2}
The family \ $(\cE_T)_{T\in\RR_{++}}$ \ of statistical experiments given in
 \eqref{cET} is LAQ at \ $-\frac{\pi^2}{2}$ \ with scaling
 \ $r_{-\frac{\pi^2}{2},T} = \frac{1}{T}$, \ $T \in \RR_{++}$, \ and with
 \begin{gather*}
  \Delta_{-\frac{\pi^2}{2}}
  = \frac{16 \int_0^1 (\cW_1(s) \, \dd \cW_2(s) - \cW_2(s) \, \dd \cW_1(s))
          - 4 \pi
            \int_0^1 (\cW_1(s) \, \dd \cW_1(s) + \cW_2(s) \, \dd \cW_2(s))}
         {\pi(\pi^2+16)} , \\
  J_{-\frac{\pi^2}{2}}
  = \frac{16}{\pi^2(\pi^2+16)}
    \int_0^1 (\cW_1(t)^2 + \cW_2(t)^2) \, \dd t,
 \end{gather*}
 where \ $(\cW_1(t), \cW_2(t))_{t\in[0,1]}$ \ is a 2-dimensional standard
 Wiener process.
\end{Pro}

\begin{Pro}\label{LAMN}
If \ $a \in (0, \infty)$ \ then the family \ $(\cE_T)_{T\in\RR_{++}}$ \ of
 statistical experiments given in \eqref{cET} is LAMN at \ $a$ \ with
 scaling \ $r_{a,T} = \ee^{-v_0(a)T}$, \ $T \in \RR_{++}$, \ and with
 \[
   \Delta_a = Z \sqrt{J_a} , \qquad
   J_a
   = \frac{(1-\ee^{-v_0(a)})^2}{2v_0(a)(v_0(a)^2+2v_0(a)-a)^2} (U^{(a)})^2 ,
 \]
 with
 \[
   U^{(a)}
   = X_0(0) + a \int_{-1}^0 \int_u^0 \ee^{-v_0(a)(s-u)} X_0(s) \, \dd s \, \dd u
     + \int_0^\infty \ee^{-v_0(a)s} \, \dd W(s) ,
 \]
 and \ $Z$ \ is a standard normally distributed random variable independent
 of \ $J_a$.
\end{Pro}

\begin{Pro}\label{PLAMN}
If \ $a \in \bigl(-\infty, -\frac{\pi^2}{2}\bigr)$ \ then the family
 \ $(\cE_T)_{T\in\RR_{++}}$ \ of statistical experiments given in \eqref{cET}
 is PLAMN at \ $a$ \ with period \ $D = \frac{\pi}{\kappa_0(a)}$, \ with
 scaling \ $r_{a,T} = \ee^{-v_0(a)T}$, \ $T \in \RR_{++}$, \ and with
 \[
   \Delta_a(d) = Z \sqrt{J_a(d)} , \qquad
   J_a(d) = \int_0^\infty \ee^{-2v_0(a)s} (V^{(a)}(d - s))^2 \, \dd s ,
   \qquad d \in \Bigl[0, \frac{\pi}{\kappa_0(a)}\Bigr) ,
 \]
 where
 \begin{align*}
  V^{(a)}(t)
  &:= X_0(0) \varphi_a(t)
      + a \int_{-1}^0 \int_u^0 
          \varphi_a(t + u - s) \ee^{-v_0(a)(s-u)} X_0(s) \, \dd s \, \dd u \\
  &\quad\:
      + \int_0^\infty \varphi_a(t - s) \ee^{-v_0(a)s} \, \dd W(s) , \qquad
  t \in \RR_+ ,
 \end{align*}
 with
 \[
   \varphi_a(t) := A_0(a) \cos(\kappa_0(a) t) + B_0(a) \sin(\kappa_0(a) t) ,
   \qquad t \in \RR ,
 \]
 and \ $Z$ \ is a standard normally distributed random variable independent
 of \ $J_a(d)$.
\end{Pro}
 
\begin{Rem}
If LAN property holds then one can construct asymptotically optimal tests,
 see, e.g., Theorem 15.4 and Addendum 15.5 of van der Vaart \cite{Vaart}.
\end{Rem}

\section{Maximum likelihood estimates}

For fixed \ $T \in \RR_{++}$, \ maximizing
 \ $\log \frac{\dd \PP_{a,T}}{\dd \PP_{0,T}}(X^{(a)}|_{[-1,T]})$ \ in
 \ $a \in \RR$ \ gives the MLE of \ $a$\ based on the observations
 \ $(X(t))_{t\in[-1,T]}$ \ having the form
 \[
   \ha_T
   = \frac{\int_0^T \int_{-1}^0 X^{(a)}(t+u) \, \dd u \, \dd X^{(a)}(t)}
          {\int_0^T \left(\int_{-1}^0 X^{(a)}(t+u) \, \dd u\right)^2 \dd t} ,
 \]
 provided that
 \ $\int_0^T \left(\int_{-1}^0 X^{(a)}(t+u) \, \dd u\right)^2 \dd t > 0$.
\ Using the SDDE \eqref{SDDE}, one can check that
 \[
   \ha_T - a
   = \frac{\int_0^T \int_{-1}^0 X^{(a)}(t+u) \, \dd u \, \dd W(t)}
          {\int_0^T \left(\int_{-1}^0 X^{(a)}(t+u) \, \dd u\right)^2 \dd t} ,
 \]
 hence
 \[
   r_{a,T}^{-1} (\ha_T - a) = \frac{\Delta_{a,T}}{J_{a,T}} .
 \]
Using the results of Section \ref{section_LALR} and the continuous mapping theorem,
 we obtain the following result.

\begin{Pro}\label{MLE}
If \ $a \in \bigl(-\frac{\pi^2}{2}, 0\bigr)$ \ then
 \[
   \sqrt{T} \, (\ha_T - a) \distr \cN(0, J_a^{-1}) \qquad
   \text{as \ $T \to \infty$,}
 \]
 where \ $J_a$ \ is given in Proposition \ref{LAN}.
 
If \ $a = 0$ \ then
 \[
   T \, (\ha_T - a) = T \, \ha_T
   \distr \frac{\int_0^1 \cW(t) \, \dd \cW(t)}{\int_0^1 \cW(t)^2 \, \dd t}
   \qquad \text{as \ $T \to \infty$,}
 \]
 where \ $(\cW(t))_{t\in[0,1]}$ \ is a standard Wiener process. 
 
If \ $a = -\frac{\pi^2}{2}$ \ then
 \begin{align*}
  &T \, (\ha_T - a) = T \, \biggl(\ha_T + \frac{\pi^2}{2}\biggr) \\
  &\distr
   \frac{16 \pi \int_0^1 (\cW_1(t) \, \dd \cW_2(t) - \cW_2(t) \, \dd \cW_1(t))
         - 4 \pi^2
           \int_0^1 (\cW_1(t) \, \dd \cW_1(t) + \cW_2(t) \, \dd \cW_2(t))}
        {16 \int_0^1 (\cW_1(t)^2 + \cW_2(t)^2) \, \dd t}
 \end{align*}
 as \ $T \to \infty$,
 \ where \ $(\cW_1(t), \cW_2(t))_{t\in[0,1]}$ \ is a
 2-dimensional standard Wiener process.
 
If \ $a \in (0, \infty)$ \ then 
 \[
   \ee^{v_0(a)T} \, (\ha_T - a) \distr \frac{Z}{\sqrt{J_a}}
   \qquad \text{as \ $T \to \infty$,}
 \]
 where \ $J_a$ \ is given in Proposition \ref{LAMN}, and \ $Z$ \ is a standard
 normally distributed random variable independent of \ $J_a$.

If \ $a \in \bigl(-\infty, -\frac{\pi^2}{2}\bigr)$ \ then for each
 \ $d \in \bigl[0, \frac{\pi}{\kappa_0(a)}\bigr)$,
 \[
   \ee^{v_0(a)(k\frac{\pi}{\kappa_0(a)}+d)} \,
   (\ha_{k\frac{\pi}{\kappa_0(a)}+d} - a)
   \distr \frac{Z}{\sqrt{J_a(d)}}
   \qquad \text{as \ $k \to \infty$,}
 \]
 where \ $J_a(d)$ \ is given in Proposition \ref{PLAMN}, and \ $Z$ \ is a
 standard normally distributed random variable independent of \ $J_a(d)$.
\end{Pro}

If LAMN property holds then we have local asymptotic minimax bound for
 arbitrary estimators, see, e.g., Le Cam and Yang \cite[6.6, Theorem 1]{LeCamYang}.
Maximum likelihood estimators attain this bound for bounded loss function, see, e.g.,
 Le Cam and Yang \cite[6.6, Remark 11]{LeCamYang}. 
Moreover, maximum likelihood estimators are asymptotically efficient in H\'ajek's
 convolution theorem sense (see, for example, Le Cam and Yang
 \cite[6.6, Theorem 3 and Remark 13]{LeCamYang} or Jeganathan \cite{Jeg}).

\section{Proofs}
\label{section_proofs}

For each \ $a \in \RR$ \ and each deterministic continuous function
 \ $(y(t))_{t\in\RR_+}$, \ consider a continuous stochastic process
 \ $(Y^{(a)}(t))_{t\in\RR_+}$ \ given by
 \begin{equation}\label{Y}
  Y^{(a)}(t)
  := y(t) X_0(0)
     + a \int_{-1}^0 \int_u^0 y(t + u - s) X_0(s) \, \dd s \, \dd u
     + \int_0^t y(t - s) \, \dd W(s)
 \end{equation}
 for \ $t \in [1, \infty)$.

\begin{Lem}\label{CS}
Let \ $(y(t))_{t\in\RR_+}$ \ be a deterministic continuous function.
Put
 \[
   Z(t) := \int_{-1}^0 \int_u^0 y(t + u - s) X_0(s) \, \dd s \, \dd u ,
   \qquad t \in [1, \infty) .
 \]
Then for each \ $T \in [1, \infty)$,
 \[
   \int_1^T Z(t)^2 \, \dd t
   \leq \int_{-1}^0 X_0(s)^2 \, \dd s \int_0^T y(v)^2 \, \dd v .
 \]
\end{Lem}

\noindent{\bf Proof.}
For each \ $t \in [1, \infty)$, \ by Fubini's theorem, 
 \[
   Z(t)
   = \int_{-1}^0 X_0(s) \int_{-1}^s y(t + u - s) \, \dd u \, \dd s
   = \int_{-1}^0 X_0(s) \int_{t-s-1}^t y(v) \, \dd v \, \dd s .
 \]
By the Cauchy--Schwarz inequality,
 \[
   Z(t)^2
   \leq \int_{-1}^0 X_0(s)^2 \, \dd s
        \int_{-1}^0 \left(\int_{t-s-1}^t y(v) \, \dd v\right)^2 \dd s .
 \]
Consequently,
 \[
   \int_1^T Z(t)^2 \, \dd t
   \leq \int_{-1}^0 X_0(s)^2 \, \dd s
        \int_1^T \int_{-1}^0
         \left(\int_{t-s-1}^t y(v) \, \dd v\right)^2 \dd s \, \dd t ,   
 \]
 where 
 \[
   \int_1^T \int_{-1}^0
    \left(\int_{t-s-1}^t y(v) \, \dd v\right)^2 \dd s \, \dd t
   = \int_{-1}^0 \int_1^T
      \left(\int_{t-s-1}^t y(v) \, \dd v\right)^2 \dd t \, \dd s .
 \]
Here
 \begin{align*}
  \int_1^T \left(\int_{t-s-1}^t y(v) \, \dd v\right)^2 \dd t
  &\leq \int_1^T \int_{t-s-1}^t y(v)^2 \, \dd v \, \dd t
   = \int_{-s}^T y(v)^2 \int_v^{v+s+1} \, \dd t \, \dd v \\
  &\leq \int_{-s}^T y(v)^2 \, \dd v
   \leq \int_0^T y(v)^2 \, \dd v 
 \end{align*}
 for all \ $s \in [-1, 0]$, \ hence we obtain the statement.
\proofend

\begin{Lem}\label{erg}
Let \ $(y(t))_{t\in\RR_+}$ \ be a deterministic continuous function with
 \ $\int_0^\infty y(t)^2 \, \dd t < \infty$.
\ Then for each \ $a \in \RR$,
 \begin{gather*}
  \frac{1}{T} \int_0^T Y^{(a)}(t) \, \dd t \stoch 0 \qquad
  \text{as \ $T \to \infty$,} \\
  \frac{1}{T} \int_0^T Y^{(a)}(t)^2 \, \dd t
  \stoch \int_0^\infty y(t)^2 \, \dd t \qquad
  \text{as \ $T \to \infty$.}
 \end{gather*}
\end{Lem}

\noindent{\bf Proof.}
Applying Lemma 4.3 of Gushchin and K\"uchler \cite{GusKuc1999} for the special
 case \ $X_0(s) = 0$, \ $s \in [-1, 0]$, \ we obtain
 \begin{gather*}
  \frac{1}{T} \int_0^T \int_0^t y(t - s) \, \dd W(s) \, \dd t \stoch 0 \qquad
  \text{as \ $T \to \infty$,} \\
  \frac{1}{T} \int_0^T \left(\int_0^t y(t - s) \, \dd W(s)\right)^2 \dd t
  \stoch \int_0^\infty y(t)^2 \, \dd t \qquad
  \text{as \ $T \to \infty$.}
 \end{gather*}
We have
 \[
   \frac{1}{T} \int_0^T Y^{(a)}(t) \, \dd t
   = \frac{1}{T} \int_0^1 Y^{(a)}(t) \, \dd t
     + X_0(0) I_1(T) + \frac{a}{T} \int_1^T Z(t) \, \dd t
     + \frac{1}{T} \int_1^T \int_0^t y(t - s) \, \dd W(s) \, \dd t
 \]
 for \ $T \in \RR_+$, \ where \ $(Z(t))_{t\in\RR_+}$ \ is given in Lemma \ref{CS}, and
 \[
   I_1(T) := \frac{1}{T} \int_1^T y(t) \, \dd t , \qquad T \in \RR_+ .
 \]
By Lemma \ref{CS},
 \begin{gather*}
  |I_1(T)| 
  \leq \sqrt{\int_0^T \frac{1}{T^2} \, \dd t
             \int_1^T y(t)^2 \, \dd t}
  = \sqrt{\frac{1}{T} \int_0^\infty y(t)^2 \, \dd t}
  \to 0 \qquad \text{as \ $T \to \infty$,} \\
  \left|\frac{1}{T} \int_1^T Z(t) \, \dd t\right|
  \leq \sqrt{\frac{1}{T} \int_1^T Z(t)^2 \, \dd t}
  \leq \sqrt{\frac{1}{T} \int_{-1}^0 X_0(s)^2 \, \dd s
             \int_0^\infty y(v)^2 \, \dd v}
  \as 0 \qquad \text{as \ $T \to \infty$,}
 \end{gather*}
 hence we obtain the first statement.
Moreover,
 \[
   \frac{1}{T} \int_0^T Y^{(a)}(t)^2 \, \dd t
   = \frac{1}{T} \int_0^1 Y^{(a)}(t)^2 \, \dd t
     + I_2(T) + 2 I_3(T)
     + \frac{1}{T}
       \int_1^T \left(\int_0^t y(t - s) \, \dd W(s)\right)^2 \dd t
 \]
 for \ $T \in \RR_+$, \ where
 \begin{gather*}
  I_2(T) := \frac{1}{T} \int_1^T (y(t) X_0(0) + a Z(t))^2 \, \dd t , \qquad
  T \in \RR_+ , \\
  I_3(T) := \frac{1}{T}
             \int_1^T (y(t) X_0(0) + a Z(t))
              \left(\int_0^t y(t - s) \, \dd W(s)\right) \dd t, \qquad
  T \in \RR_+ .
 \end{gather*}
Again by Lemma \ref{CS},
 \begin{align*}
  0 \leq I_2(T)
  &\leq \frac{1}{T} \int_1^T 2 (y(t)^2 X_0(0)^2 + a^2 Z(t)^2) \, \dd t \\
  &\leq \frac{2X_0(0)^2}{T} \int_0^\infty y(t)^2 \, \dd t
        + \frac{2a^2}{T} \int_{-1}^0 X_0(s)^2 \, \dd s
          \int_0^\infty y(v)^2 \, \dd v
  \as 0 \qquad \text{as \ $T \to \infty$,}
 \end{align*}
 and
 \begin{align*}
  |I_3(T)| 
  &\leq \frac{2}{T} 
        \sqrt{\int_1^T (y(t)^2 X_0(0)^2 + a^2 Z(t)^2)^2 \, \dd t 
              \int_1^T \left(\int_0^t y(t - s) \, \dd W(s)\right)^2 \dd t} \\
  &= 2 \sqrt{\frac{I_2(T)}{T} 
             \int_1^T \left(\int_0^t y(t - s) \, \dd W(s)\right)^2 \dd t}
  \stoch 0 \qquad \text{as \ $T \to \infty$,}
 \end{align*}
 hence we obtain the second statement.
\proofend

\begin{Lem}\label{super1}
Let \ $w \in \RR_{++}$ \ and \ $y(t) := \ee^{wt}$, \ $t \in \RR_+$.
\ Then for each \ $a \in \RR$,
 \begin{gather*}
  \ee^{-wt} Y^{(a)}(t) \as U^{(a)}_w , \qquad \text{as \ $t \to \infty$,} \\
  \ee^{-2wT} \int_0^T (Y^{(a)}(t))^2 \, \dd t \as \frac{1}{2w} (U^{(a)}_w)^2 ,
  \qquad \text{as \ $T \to \infty$,}
 \end{gather*}
 with
 \[
   U^{(a)}_w 
   := X_0(0)
      + a \int_{-1}^0 \int_u^0 \ee^{w(u-s)} X_0(s) \, \dd s \, \dd u 
      + \int_0^\infty \ee^{-ws} \, \dd W(s) .
 \]
\end{Lem}

\noindent{\bf Proof.}
For each \ $t \in [1, \infty)$, \ we have
 \[
   \ee^{-wt} Y^{(a)}(t)
   = X_0(0)
      + a \int_{-1}^0 \int_u^0 \ee^{w(u-s)} X_0(s) \, \dd s \, \dd u 
      + \int_0^t \ee^{-ws} \, \dd W(s) ,
 \]
 hence we obtain the first convergence.
The second convergence follows by L'H\^ospital's rule.
\proofend

\begin{Lem}\label{super2}
Let \ $w \in \RR_{++}$, \ $\kappa \in \RR$, \ and
 \ $y(t) := \varphi(t) \ee^{wt}$, \ $t \in \RR_+$, \ with
 \ $\varphi(t) = \cos(\kappa t)$, \ $t \in \RR_+$, \ or
 \ $\varphi(t) = \sin(\kappa t)$, \ $t \in \RR_+$.
\ Then for each \ $a \in \RR$,
 \begin{gather*}
  \ee^{-wt} Y^{(a)}(t) - V^{(a)}_w(t) \as 0 , \qquad
  \text{as \ $t \to \infty$,} \\
  \ee^{-2wT} \int_0^T (Y^{(a)}(t))^2 \, \dd t
  - \int_0^\infty \ee^{-2wt} (V^{(a)}_w(T-t))^2 \, \dd t
  \stoch 0 , \qquad \text{as \ $T \to \infty$,}
 \end{gather*}
 with
 \[
   V^{(a)}_w(t)
   := X_0(0) \varphi(t)
      + a \int_{-1}^0 \int_u^0 
           \varphi(t + u - s) \ee^{w(u-s)} X_0(s) \, \dd s \, \dd u 
      + \int_0^\infty \varphi(t - s) \ee^{-ws} \, \dd W(s)
 \]
 for \ $t \in \RR$.
\end{Lem}

\noindent{\bf Proof.}
Note that for each \ $t \in [1, \infty)$,
 \begin{equation}\label{Vaw}
   \ee^{-wt} Y^{(a)}(t) - V^{(a)}_w(t)
   = - \int_t^\infty \varphi(t - s) \ee^{-ws} \, \dd W(s) ,
 \end{equation}
 which obviously tends almost surely to zero as \ $t \to \infty$, \ hence we
 obtain the first convergence.
 
In order to prove the second convergence, observe that for each
 \ $T \in \RR_+$,
 \[
   \int_0^\infty \ee^{-2wt} (V^{(a)}_w(T-t))^2 \, \dd t
   = \int_0^T \ee^{-2wt} (V^{(a)}_w(T-t))^2 \, \dd t
     + \int_T^\infty \ee^{-2wt} (V^{(a)}_w(T-t))^2 \, \dd t ,
 \]
 where
 \[
   \int_0^T \ee^{-2wt} (V^{(a)}_w(T-t))^2 \, \dd t
   = \int_0^T \ee^{-2w(T-t)} (V^{(a)}_w(t))^2 \, \dd t
   = \ee^{-2wT} \int_0^T (\ee^{wt} V^{(a)}_w(t))^2 \, \dd t ,
 \]
 hence
 \begin{align*}
  &\ee^{-2wT} \int_0^T (Y^{(a)}(t))^2 \, \dd t
   - \int_0^\infty \ee^{-2wt} (V^{(a)}_w(T-t))^2 \, \dd t \\
  &= \ee^{-2wT}
     \int_0^T
      \bigl[ (Y^{(a)}(t))^2 - (\ee^{wt} V^{(a)}_w(t))^2 \bigr] \dd t
     - \int_T^\infty \ee^{-2wt} (V^{(a)}_w(T-t))^2 \, \dd t \\
  &= I_0(T) + I_1(T) + 2 I_2(T) - I_3(T) 
 \end{align*}
 with
 \begin{align*}
  I_0(T)
  &:= \ee^{-2wT} \int_0^1 (Y^{(a)}(t) - \ee^{wt} V^{(a)}_w(t))^2 \, \dd t , \\
  I_1(T)
  &:= \ee^{-2wT} \int_1^T (Y^{(a)}(t) - \ee^{wt} V^{(a)}_w(t))^2 \, \dd t , \\
  I_2(T)
  &:= \ee^{-2wT}
      \int_0^T
       (Y^{(a)}(t) - \ee^{wt} V^{(a)}_w(t)) \ee^{wt} V^{(a)}_w(t) \, \dd t , \\
  I_3(T)
  &:= \int_T^\infty \ee^{-2wt} (V^{(a)}_w(T-t))^2 \, \dd t .
 \end{align*}
The processes \ $(Y^{(a)}(t))_{t\in\RR_+}$ \ and \ $(V^{(a)}_w(t))_{t\in\RR_+}$
 \ are continuous, hence
 \[
   \EE(|I_0(T)|)
   = \ee^{-2wT} \int_0^1 \EE\bigl[(Y^{(a)}(t) - \ee^{wt} V^{(a)}_w(t))^2\bigr] \dd t
   \to 0 \qquad \text{as \ $T \to \infty$,}
 \]
 implying \ $I_0(T) \stoch 0$ \ as \ $T \to \infty$.
\ By \eqref{Vaw},
 \[
   I_1(T)
   = \ee^{-2wT}
     \int_1^T
      \ee^{2wt}
      \left(\int_t^\infty \varphi(t - s) \ee^{-ws} \, \dd W(s)\right)^2
      \dd t ,
 \]
 hence
 \begin{align*}
  \EE(|I_1(T)|)
  &= \ee^{-2wT}
     \int_1^T
      \ee^{2wt}
      \int_t^\infty \varphi(t - s)^2 \ee^{-2ws} \, \dd s \, \dd t
   \leq \ee^{-2wT}
        \int_1^T \ee^{2wt} \int_t^\infty \ee^{-2ws} \, \dd s \, \dd t \\
  &= \frac{T}{2w} \ee^{-2wT} \to 0 \qquad \text{as \ $T \to \infty$,}
 \end{align*}
 implying \ $I_1(T) \stoch 0$ \ as \ $T \to \infty$.
\ Moreover, by the Cauchy--Schwarz inequality,
 \[
   |I_2(T)|
   \leq \sqrt{I_1(T) \ee^{-2wT}
              \int_0^T (\ee^{wt} V^{(a)}_w(t))^2 \, \dd t}
 \]
 with
 \[
   \ee^{-2wT} \int_0^T (\ee^{wt} V^{(a)}_w(t))^2 \, \dd t
   = \int_0^T \ee^{-2wt} V^{(a)}_w(T - t)^2 \, \dd t
   \leq \frac{1}{2w} \sup_{t\in\RR} (V^{(a)}_w(t))^2 ,
 \]
 where \ $\sup_{t\in\RR} (V^{(a)}_w(t))^2 < \infty$ \ almost surely, since
 \ $(V^{(a)}_w(t))_{t\in\RR}$ \ is a continuous and periodic process.
Consequently, \ $I_2(T) \stoch 0$ \ as \ $T \to \infty$.
\ Finally,
 \[
   |I_3(T)| \leq \frac{\ee^{-2wT}}{2w} \sup_{t\in\RR} (V^{(a)}_w(t))^2
   \as 0 \qquad \text{as \ $T \to \infty$,}
 \]
 hence we obtain the second convergence of the statement.
\proofend

\noindent{\bf Proof of Proposition \ref{LAN}.}
For each \ $t \in [1, \infty)$, \ by \eqref{solution}, we have
 \begin{align*}
  \int_{-1}^0 X^{(a)}(t + u) \, \dd u
  &= X_0(0) \int_{-1}^0 x_{0,a}(t + u) \, \dd u
     + a \int_{-1}^0 \int_{-1}^0 \int_v^0
          x_{0,a}(t + u + v - s) X_0(s) \, \dd s \, \dd v \, \dd u \\
  &\quad
     + \int_{-1}^0 \int_0^{t+u} x_{0,a}(t + u - s) \, \dd W(s) \, \dd u .
 \end{align*}
Here we have
 \begin{align*}
  \int_{-1}^0 \int_{-1}^0 \int_v^0
   x_{0,a}(t + u + v - s) X_0(s) \, \dd s \, \dd v \, \dd u
  &= \int_{-1}^0 \int_{-1}^0 \int_v^0
      x_{0,a}(t + u + v - s) X_0(s) \, \dd s \, \dd u \, \dd v \\
  &= \int_{-1}^0 \int_v^0 X_0(s) \int_{-1}^0
      x_{0,a}(t + u + v - s) \, \dd u \, \dd s \, \dd v ,
 \end{align*}
 and
 \begin{align*}
  &\int_{-1}^0 \int_0^{t+u} x_{0,a}(t + u - s) \, \dd W(s) \, \dd u \\
  &= \int_0^{t-1} \int_{-1}^0 x_{0,a}(t + u - s) \, \dd u \, \dd W(s)
     + \int_{t-1}^t \int_{s-t}^0 x_{0,a}(t + u - s) \, \dd u \, \dd W(s) \\
  &= \int_0^t \int_{-1}^0 x_{0,a}(t + u - s) \, \dd u \, \dd W(s) ,
 \end{align*}
 since \ $t \in [1, \infty)$, \ $s \in [t-1, t]$ \ and \ $u \in [-1, s-t)$
 \ imply \ $t + u - s \in [-1, 0)$, \ and hence \ $x_{0,a}(t + u - s) = 0$.
\ Consequently, the process
 \ $\bigl(\int_{-1}^0 X^{(a)}(t + u) \, \dd u\bigr)_{t\in\RR_+}$
 \ has a representation \eqref{Y} with
 \[
   y(t) = \int_{-1}^0 x_{0,a}(t + u) \, \dd u , \qquad t \in \RR_+ .
 \]
Assumption \ $a \in \bigl(- \frac{\pi^2}{2}, 0\bigr)$ \ implies
 \ $v_0(a) < 0$, \ and hence
 \ $\int_0^\infty x_{0,a}(t)^2 \, \dd t < \infty$ \ holds.
Thus
 \begin{align*}
  \int_1^\infty y(t)^2 \, \dd t
  &= \int_1^\infty
        \left(\int_{-1}^0 x_{0,a}(t + u) \, \dd u\right)^2 \dd t \\
  &= \int_{-1}^0 \int_{-1}^0 \int_1^\infty
        x_{0,a}(t + u) x_{0,a}(t + v) \, \dd t \, \dd u \, \dd v
   \leq \int_0^\infty x_{0,a}(t)^2 \, \dd t ,
 \end{align*}
 since
 \begin{align*}
  &\left|\int_1^\infty x_{0,a}(t + u) x_{0,a}(t + v) \, \dd t\right|
   \leq \sqrt{\int_1^\infty x_{0,a}(t + u)^2 \, \dd t
              \int_1^\infty x_{0,a}(t + v)^2 \, \dd t} \\
  &\qquad\qquad\qquad\qquad\qquad
   = \sqrt{\int_{1+u}^\infty x_{0,a}(s + u)^2 \, \dd s
           \int_{1+v}^\infty x_{0,a}(s + v)^2 \, \dd s}
   \leq \int_0^\infty x_{0,a}(t)^2 \, \dd t .
 \end{align*}
Consequently,
 \ $\int_0^\infty y(t)^2 \, \dd t
    \leq \int_0^1 y(t)^2 \, \dd t + \int_0^\infty x_{0,a}(t)^2 \, \dd t < \infty$,
 \ thus we can apply Lemma \ref{erg} to obtain
 \[
   J_{a,T}
   = \frac{1}{T}
     \int_0^T \left(\int_{-1}^0 X^{(a)}(t + u) \, \dd u \right)^2 \dd t
   \stoch
   \int_0^\infty
    \left(\int_{-1}^0 x_{0,a}(t + u) \, \dd u\right)^2 \dd t
   = J_a
 \]
 as \ $T \to \infty$.
\ Moreover, the process 
 \[
   M^{(a)}(T) := \int_0^T \int_{-1}^0 X^{(a)}(t+u) \, \dd u \, \dd W(t) ,
   \qquad T \in \RR_+ ,
 \]
 is a continuous martingale with \ $M^{(a)}(0) = 0$ \ and with quadratic
 variation
 \[
   \langle M^{(a)} \rangle(T)
   = \int_0^T \left(\int_{-1}^0 X^{(a)}(t+u) \, \dd u\right)^2 \dd t ,
 \]
 hence, Theorem VIII.5.42 of Jacod and Shiryaev \cite{JSh} yields the statement.
\proofend

\noindent{\bf Proof of Proposition \ref{LAQ1}.}
We have
 \[
   \Delta_{0,T}
   = \frac{1}{T} \int_0^T \int_{-1}^0 X^{(0)}(t + u) \, \dd u \, \dd W(t) ,
   \qquad T \in \RR_{++} .
 \]
As in the proof of Proposition \ref{LAN}, for each \ $t \in [1, \infty)$,
 \ we obtain
 \[
   \int_{-1}^0 X^{(0)}(t + u) \, \dd u
   = X_0(0) \int_{-1}^0 x_{0,0}(t + u) \, \dd u
     + \int_0^t \int_{-1}^0 x_{0,0}(t + u - s) \, \dd u \, \dd W(s) .
 \]
Here we have
 \[
   \int_{-1}^0 x_{0,0}(t + u) \, \dd u = 1 , \qquad
   \int_{-1}^0 x_{0,0}(t + u - s) \, \dd u
   = \begin{cases}
      1 , & \text{for \ $s \in [0, t - 1]$,} \\
      t - s , & \text{for \ $s \in [t - 1, t]$,}
     \end{cases}
 \]
 hence
 \begin{align*}
  \int_{-1}^0 X^{(0)}(t + u) \, \dd u
  &= X_0(0) + \int_0^{t-1} \dd W(s) + \int_{t-1}^t (t - s) \, \dd W(s) \\
  &= X_0(0) + W(t) + \int_{t-1}^t (t - s - 1) \, \dd W(s)
   = W(t) + \overline{X}(t) ,
 \end{align*}
 where \ $\EE\bigl(T^{-2} \int_0^T \overline{X}(t)^2 \, \dd t\bigr) \to 0$
 \ as \ $T \to \infty$.
\ For each \ $T \in \RR_{++}$, \ consider the process
 \[
   W^T(s) := \frac{1}{\sqrt{T}} W(Ts) , \qquad s \in [0, 1] .
 \]
Then we have
 \begin{align*}
  \Delta_{0,T}
  &= \int_0^1 W^T(t) \, \dd W^T(t)
     + \frac{1}{T} \int_0^T \overline{X}(t) \, \dd W(t) , \\
  J_{0,T}
  &= \int_0^1 W^T(t)^2 \, \dd t
     + \frac{2}{T^2} \int_0^T W(t) \overline{X}(t) \, \dd t
     + \frac{1}{T^2} \int_0^T \overline{X}(t)^2 \, \dd t .
 \end{align*}
Here
 \[
   \frac{1}{T} \int_0^T \overline{X}(t) \, \dd W(t) \stoch 0 , \qquad
   \frac{1}{T^2} \int_0^T \overline{X}(t)^2 \, \dd t \stoch 0 
 \]
 as \ $T \to \infty$, \ since
 \[
   \EE\left[\left(\frac{1}{T}
                  \int_0^T \overline{X}(t) \, \dd W(t)\right)^2\right]
   = \frac{1}{T^2} \int_0^T \EE(\overline{X}(t)^2) \, \dd t \to 0 .
 \]
By the functional central limit theorem,
 \[
   W^T \distr \cW \qquad \text{as \ $T \to \infty$,}
 \]
 hence
 \begin{align*}
  \left|\frac{1}{T^2} \int_0^T W(t) \overline{X}(t) \, \dd t\right|
  &\leq \sqrt{\left(\frac{1}{T^2} \int_0^T W(t)^2 \, \dd t\right)
              \left(\frac{1}{T^2}
                    \int_0^T \overline{X}(t)^2 \, \dd t\right)} \\
  &= \sqrt{\left(\int_0^1 W^T(t)^2 \, \dd t\right)
           \left(\frac{1}{T^2} \int_0^T \overline{X}(t)^2 \, \dd t\right)}
   \stoch 0 \qquad \text{as \ $T \to \infty$,}
 \end{align*}
 and the claim follows from Corollary 4.12 in Gushchin and K\"uchler
 \cite{GusKuc1999}.
\proofend

\noindent{\bf Proof of Proposition \ref{LAQ2}.}
We have
 \begin{gather*}
  \Delta_{-\frac{\pi^2}{2},T}
  = \frac{1}{T}
    \int_0^T \int_{-1}^0 X^{(-\pi^2/2)}(t + u) \, \dd u \, \dd W(t) ,
  \qquad T \in \RR_{++} , \\
  J_{-\frac{\pi^2}{2},T}
  = \frac{1}{T^2}
    \int_0^T
     \biggl(\int_{-1}^0 X^{(-\pi^2/2)}(t + u) \, \dd u\biggr)^2
     \dd t ,
  \qquad T \in \RR_{++} .
 \end{gather*}
As in the proof of Proposition \ref{LAN}, for each \ $t \in [1, \infty)$,
 \ we have
 \begin{align*}
  \int_{-1}^0 X^{(-\pi^2/2)}(t + u) \, \dd u
  &= X_0(0) \int_{-1}^0 x_{0,-\frac{\pi^2}{2}}(t + u) \, \dd u
     + \int_0^t \int_{-1}^0
        x_{0,-\frac{\pi^2}{2}}(t + u - s) \, \dd u \, \dd W(s) \\
  &\quad 
     - \frac{\pi^2}{2}
       \int_{-1}^0 \int_v^0 X_0(s) \int_{-1}^0
        x_{0,-\frac{\pi^2}{2}}(t + u + v - s) \, \dd u \, \dd s \, \dd v .
 \end{align*}
We have \ $v_0\bigl(-\frac{\pi^2}{2}\bigr) = 0$ \ and
 \ $\kappa_0\bigl(-\frac{\pi^2}{2}\bigr) = \pi$, \ hence
 \ $A_0\bigl(-\frac{\pi^2}{2}\bigr) = \frac{16}{\pi^2+16}$ \ and
 \ $B_0\bigl(-\frac{\pi^2}{2}\bigr) = \frac{4\pi}{\pi^2+16}$.
\ Consequently, by Lemma \ref{res}, there exists
 \ $\gamma \in (-\infty, 0)$ \ such that
 \[
   x_{0,-\frac{\pi^2}{2}}(t)
   = \frac{16\cos(\pi t)+4\pi\sin(\pi t)}{\pi^2+16}
     + \oo(\ee^{\gamma t}) ,
   \qquad \text{as \ $t \to \infty$,}
 \]
 and hence
 \begin{align*}
  \int_{-1}^0 X^{(-\pi^2/2)}(t + u) \, \dd u
  &= \int_0^t \int_{-1}^0 
      \frac{16\cos(\pi(t+u-s))+4\pi\sin(\pi(t+u-s))}{\pi^2+16}
      \, \dd u \, \dd W(s)
     + \overline{X}(t) \\
  &= \int_0^t
      \frac{32 \sin(\pi(t-s)) - 8 \pi \cos(\pi(t-s))}{\pi(\pi^2+16)}
      \, \dd W(s)
     + \overline{X}(t) ,
 \end{align*}
 where
 \ $T^{-2} \int_0^T \overline{X}(t)^2 \, \dd t \stoch 0$ \ as
 \ $T \to \infty$.
\ Introducing
 \[
   X_1(t) := \int_0^t \cos(\pi s) \, \dd W(s) , \qquad
   X_2(t) := \int_0^t \sin(\pi s) \, \dd W(s) , \qquad t \in \RR_+ ,
 \]
 we obtain
 \begin{align*}
  &\int_{-1}^0 X^{(-\pi^2/2)}(t + u) \, \dd u \\
  &\qquad\qquad
   = \frac{32 X_1(t) \sin(\pi t) - 32 X_2(t) \cos(\pi t)
           - 8 \pi X_1(t) \cos(\pi t) - 8 \pi X_2(t) \sin(\pi t)}
          {\pi(\pi^2+16)}
     + \overline{X}(t) .
 \end{align*}
For each \ $T \in \RR_{++}$, \ consider the following processes on \ $[0, 1]$:
 \begin{gather*}
  W^T(s) := \frac{1}{\sqrt{T}} W(Ts) , \\
  X_1^T(s) := \frac{1}{\sqrt{T}} X_1(Ts)
            = \int_0^s \cos(\pi T s) \, \dd W^T(s) , \\
  X_2^T(s) := \frac{1}{\sqrt{T}} X_2(Ts)
            = \int_0^s \sin(\pi T s) \, \dd W^T(s), \\
  X^T(s) := \frac{32 X_1^T(s) \sin(\pi Ts) - 32 X_2^T(s) \cos(\pi Ts)
             - 8 \pi X_1(s) \cos(\pi Ts) - 8 \pi X_2(s) \sin(\pi Ts)}
                 {\pi(\pi^2+16)} .
 \end{gather*}
Then, for each \ $T \in \RR_{++}$, \ we have
 \[
   \int_{-1}^0 X^{(-\pi^2/2)}(t + u) \, \dd u
   = \sqrt{T} X^T\Bigl(\frac{t}{T}\Bigr) + \overline{X}(t) ,
 \]
 and hence,
 \begin{align*}
  \Delta_{-\frac{\pi^2}{2},T}
  &= \frac{1}{\sqrt{T}} \int_0^T X^T\Bigl(\frac{t}{T}\Bigr) \, \dd W(t)
     + I_1(T)
   = \int_0^1 X^T(s) \, \dd W^T(s) + I_1(T) , \\
  J_{-\frac{\pi^2}{2},T}
  &= \frac{1}{T}
     \int_0^T X^T\Bigl(\frac{t}{T}\Bigr)^2 \dd t
     + 2 I_2(T) + I_3(T)
   = \int_0^1 X^T(s)^2 \, \dd s + 2 I_2(T) + I_3(T) ,
 \end{align*}
 with
 \[
   I_1(T) := \frac{1}{T} \int_0^T \overline{X}(t) \, \dd W(t) , \quad
   I_2(T) := \frac{1}{T^{3/2}}
             \int_0^T
              X^T\Bigl(\frac{t}{T}\Bigr) \overline{X}(t) \, \dd t , \quad
   I_3(T) := \frac{1}{T^2} \int_0^T \overline{X}(t)^2 \, \dd t .
 \]
Introducing the process
 \[
   Y^T(t) := \int_0^t X^T(s) \, \dd W^T(s) ,
   \qquad t \in \RR_+ , \qquad T \in \RR_{++} ,
 \]
 we have
 \[
   \int_0^t X^T(s)^2 \, \dd s = [Y^T, Y^T](t) ,
   \qquad t \in \RR_+ , \qquad T \in \RR_{++} ,
 \]
 where \ $([U, V](t))_{t\in\RR_+}$ \ denotes the quadratic covariation process
 of the processes \ $(U(t))_{t\in\RR_+}$ \ and \ $(V(t))_{t\in\RR_+}$.
\ Moreover,
 \[
   Y^T(t)
   = \frac{32 \int_0^t (X_1^T(s) \, \dd X_2^T(s) - X_2^T(s) \, \dd X_1^T(s))
           - 8 \pi
             \int_0^t (X_1^T(s) \, \dd X_1^T(s) + X_2^T(s) \, \dd X_2^T(s))}
          {\pi(\pi^2+16)}
 \]
 for \ $t \in \RR_+$.
\ By the functional central limit theorem,
 \[
   (X_1^T, X_2^T) \distr \frac{1}{\sqrt{2}} (\cW_1, \cW_2) \qquad 
   \text{as \ $T \to \infty$,}
 \]
 hence
 \[
   Y^T \distr \cY \qquad \text{as \ $T \to \infty$}
 \]
 with
 \[
   \cY(t)
   := \frac{16 \int_0^t (\cW_1(s) \, \dd \cW_2(s) - \cW_2(s) \, \dd \cW_1(s))
            - 4 \pi
              \int_0^t (\cW_1(s) \, \dd \cW_1(s) + \cW_2(s) \, \dd \cW_2(s))}
           {\pi(\pi^2+16)}
 \]
 for \ $t \in \RR_+$.
\ Further, by Corollary 4.12 in Gushchin and K\"uchler \cite{GusKuc1999},
 \[
   (Y^T(1), [Y^T, Y^T](1)) \distr (\cY(1), [\cY, \cY](1)) \qquad 
   \text{as \ $T \to \infty$.}
 \]
Here we have
 \begin{align*}
  [\cY, \cY](1)
  &= \frac{\int_0^1 (16 \cW_1(s) - 4 \pi \cW_2(s))^2 \dd s
           + \int_0^1 (16 \cW_2(s) + 4 \pi \cW_1(s))^2 \dd s} 
          {\pi^2(\pi^2+16)^2} \\
  &= \frac{16}{\pi^2(\pi^2+16)}
     \int_0^1 (\cW_1(s)^2 + \cW_2(s)^2) \, \dd s .
 \end{align*}
Recall that \ $I_3(T) \stoch 0$ \ as \ $T \to \infty$, \ which also implies
 \ $I_1(T) \stoch 0$ \ as \ $T \to \infty$.
\ Finally,
 \[
   |I_2(T)|
   \leq \sqrt{\frac{1}{T^3}
              \int_0^T X^T\Bigl(\frac{t}{T}\Bigr)^2 \dd t
              \int_0^T \overline{X}(t)^2 \, \dd t}
   = \sqrt{\frac{1}{T^2}
           \int_0^1 X^T(s)^2 \, \dd s
           \int_0^T \overline{X}(t)^2 \, \dd t}
   \stoch 0
 \]
 as \ $T \to \infty$, \ and the claim follows.
\proofend

\noindent{\bf Proof of Proposition \ref{LAMN}.}
We have
 \[
   J_{a,T}
   = \ee^{-2v_0(a)T}
     \int_0^T \left(\int_{-1}^0 X^{(a)}(t + u) \, \dd u \right)^2 \dd t
   \qquad T \in \RR_+ .
 \]
The process
 \ $\bigl(\int_{-1}^0 X^{(a)}(t + u) \, \dd u\bigr)_{t\in[1,\infty)}$
 \ has a representation \eqref{Y} with
 \ $y(t) = \int_{-1}^0 x_{0,a}(t + u) \, \dd u$, \ $t \in \RR_+$, \ see
 the proof of Proposition \ref{LAN}.
The assumption \ $a \in (0, \infty)$ \ implies \ $v_0(a) > 0$ \ and
 \ $v_1(a) < 0$, \ hence by Lemma \ref{res}, there exists
 \ $\gamma \in (v_1(a), 0)$ \ such that
 \[
   x_{0,a}(t)
   = \frac{v_0(a)}{v_0(a)^2+2v_0(a)-a} \, \ee^{v_0(a)t}
     + \oo(\ee^{\gamma t}) ,
   \qquad \text{as \ $t \to \infty$.}
 \]
Consequently,
 \[
   \int_{-1}^0 x_{0,a}(t + u) \, \dd u
   = \frac{1-\ee^{-v_0(a)}}{v_0(a)^2+2v_0(a)-a} \, \ee^{v_0(a)t}
     + \oo(\ee^{\gamma t}) ,
   \qquad \text{as \ $t \to \infty$.}
 \]
Applying Lemma \ref{super1}, we obtain
 \[
   J_{a,T}
   \stoch
   \frac{1}{2v_0(a)}
   \left(\frac{1-\ee^{-v_0(a)}}{v_0(a)^2+2v_0(a)-a}\right)^2
   (U^{(a)})^2
   = J_a \qquad \text{as \ $T \to \infty$.}
 \]
Theorem VIII.5.42 of Jacod and Shiryaev \cite{JSh} yields the statement.
\proofend

\noindent{\bf Proof of Proposition \ref{PLAMN}.}
We have again
 \[
   J_{a,T}
   = \ee^{-2v_0(a)T}
     \int_0^T \left(\int_{-1}^0 X^{(a)}(t + u) \, \dd u \right)^2 \dd t
   \qquad T \in \RR_+ ,
 \]
 and the process
 \ $\bigl(\int_{-1}^0 X^{(a)}(t + u) \, \dd u\bigr)_{t\in[1,\infty)}$
 \ has a representation \eqref{Y} with
 \ $y(t) = \int_{-1}^0 x_{0,a}(t + u) \, \dd u$, \ $t \in \RR_+$, \ see
 the proof of Proposition \ref{LAN}.
The assumption \ $a \in \bigl(-\infty, -\frac{\pi^2}{2}\bigr)$ \ implies
 \ $v_0(a) > 0$ \ and \ $v_0(a) \notin \Lambda_a$, \ hence by Lemma 
 \ref{res}, there exists \ $\gamma \in (0, v_0(a))$ \ such that
 \[
   x_{0,a}(t)
   = \varphi_a(t) \ee^{v_0(a)t}
     + \oo(\ee^{\gamma t}) ,
   \qquad \text{as \ $t \to \infty$.}
 \]
Applying Lemma \ref{super2}, we obtain
 \[
   J_{a,T} - J_a(T) \stoch 0 , \qquad \text{as \ $T \to \infty$.}
 \]
The process \ $(J_a(t))_{t\in\RR_+}$ \ is periodic with period
 \ $D = \frac{\pi}{\kappa_0(a)}$.
\proofend

\end{document}